\catcode`\@=11

\font\tenmsb=msbm10
\font\sevenmsb=msbm7 \font\fivemsb=msbm5  \newfam\msbfam
\textfont\msbfam=\tenmsb\scriptfont\msbfam=\sevenmsb
\scriptscriptfont\msbfam=\fivemsb

\def\hexnumber@#1{\ifnum#1<10 \number#1\else \ifnum#1=10 A\else\ifnum#1=11
 B\else\ifnum#1=12 C\else \ifnum#1=13 D\else\ifnum#1=14 E\else\ifnum#1=15
 F\fi\fi\fi\fi\fi\fi\fi}
 \def\msb@{\hexnumber@\msbfam}

\mathchardef\hbar="0\msb@7E

\def\Bbb{\ifmmode\let\next\Bbb@\else
\def\next{\errmessage{Use \string\Bbb\space only in math mode}}\fi\next}
\def\Bbb@#1{{\Bbb@@{#1}}} \def\Bbb@@#1{\fam\msbfam#1}

\catcode`\@=\active

\def\del{\partial}
\def\CR{\hbox{{$\cal R$}}}
\def\CM{\hbox{{$\cal M$}}}

\def\Z{{\Bbb Z}}
\def\H{{\Bbb H}}

\def\CC{{\cal C}}

\def\isom{{\cong}}
\def\eps{{\epsilon}}
\def\tens{\mathop{\otimes}}
\def\la{{\triangleright}}

\def\proof{\goodbreak\noindent{\bf Proof\quad}}
\def\endproof{{\ $\hbox{$\sqcup$}\llap{\hbox{$\sqcap$}}$}\bigskip }

\def\ev{{\rm ev}}\def\coev{{\rm coev}}
\def\End{{\rm End}}
\def\id{{\rm id}}
\def\<{\langle}
\def\>{\rangle}

\def\und#1{{\underline {#1}}}

\def\o{{}_{\scriptscriptstyle(1)}}
\def\t{{}_{\scriptscriptstyle(2)}}
\def\fo{{}_{\scriptscriptstyle(4)}}
\def\fiv{{}_{\scriptscriptstyle(5)}}
\def\th{{}_{\scriptscriptstyle(3)}}
\def\bo{{}^{\bar{\scriptscriptstyle(1)}}}
\def\bt{{}^{\bar{\scriptscriptstyle(2)}}}

\def\note#1{{}}

\def\nquad{{\!\!\!\!\!\!}}

\def\equad{\nquad}
\def\cmath#1{\[\begin{array}{c} #1 \end{array}\]}
\def\eqn#1#2{\begin{equation}#2\label{#1}\end{equation}}
\def\ceqn#1#2{\begin{equation}\label{#1}
\begin{array}{c}#2\end{array}\end{equation}}
\def\align#1{\begin{eqnarray*}#1\end{eqnarray*}}

\def\kfg{$k_FG$}

\documentstyle[11pt]{article}
\newtheorem{lemma}{Lemma}[section]
\newtheorem{propos}[lemma]{Proposition}
\newtheorem{example}[lemma]{Example}

\newtheorem{corol}[lemma]{Corollary}
\newtheorem{defin}[lemma]{Definition}

\begin{document}\baselineskip 22pt

{\ }\qquad  \hskip 4.3in DAMTP/97-138
\vspace{.2in}

\begin{center} {\LARGE QUASIALGEBRA STRUCTURE OF THE OCTONIONS}
\\ \baselineskip 13pt{\ }
{\ }\\ Helena Albuquerque\footnote{Supported by CMUC-JNICT and by
Praxis 2/2.1/Mat7458/94} \\{\ }\\ Departamento de
Matematica-Faculdade de Ciencias e Tecnologia\\ Universidade de
Coimbra, Apartado 3008\\ 3000 Coimbra, Portugal\\ +\\ Shahn
Majid\footnote{Royal Society University Research Fellow and Fellow
of Pembroke College, Cambridge}\\ {\ }\\ Department of Applied
Mathematics \& Theoretical Physics\\ University of Cambridge,
Cambridge CB3 9EW, UK\\ www.damtp.cam.ac.uk/user/majid
\end{center}
\begin{center}
December, 1997 --  February, 1998
\end{center}

\vspace{10pt}
\begin{quote}\baselineskip 13pt
\noindent{\bf Abstract}
We show that the octonions are a twisting of the group algebra of
$\Z_2\times\Z_2\times\Z_2$ in the quasitensor category of
representations of a quasi-Hopf algebra associated to a group
3-cocycle. We consider general quasi-associative algebras of this
type and some general constructions for them, including
quasi-linear algebra and representation theory, and an automorphism
quasi-Hopf algebra. Other examples include the higher $2^n$-onion
Cayley algebras and examples associated to Hadamard matrices.

\end{quote}
\baselineskip 22pt

\section{Introduction}

In this paper we provide a natural setting for the octonion
algebra, namely as an algebra in a quasitensor category. Such
categories have a tensor product and associativity isomorphisms
$V\tens (W\tens Z)\isom (V\tens W)\tens Z$ for any three objects,
but these need not, however, be the trivial vector space
isomorphisms. These categories also arise naturally as the
representation categories of quasi-Hopf algebras\cite{Dri:qua}. Our
first result is to identify the correct `octonion generating
quasi-Hopf algebra' in the category of representations of which the
octonions live. The categorical point of view then provides further
general constructions on the octonions. Moreover, the framework has
many other interesting quasi-associative algebras beyond these.

Our general construction of the octonions mirrors, for discrete
groups, Drinfeld's construction of the quantum groups $U_q(g)$.
Namely, we consider finite group function algebras $k(G)$ ($k$ a
field) regarded trivially as quasi-Hopf algebras $(k(G),\phi)$
where $\phi$ is a group 3-cocycle on $G$. However, for the
octonions, the cocycle is a coboundary and can be identified as the
result of twisting $k(G)$ by a 2-cochain $F$. In an extension of
Drinfeld's theory of twisting or `gauge equivalence'\cite{Dri:qua},
any algebra on which $k(G)$ acts also has to be twisted to remain a
module-algebra. In this sense the octonions are gauge-equivalent or
the twisting of the group algebra of $\Z_2\times\Z_2\times \Z_2$ by
a 2-cochain, which we provide. This accounts for many of the
properties of octonions as gauge-equivalent to properties of a
group algebra.

Section~2 recalls preliminaries about quasi-Hopf algebras and
quasitensor categories. As a modest result, we give in detail the
twisting theory of module algebras associated to the twisting of
quasi-Hopf algebras. We then study the case of quasi-Hopf algebras
$(k(G),\phi)$ and hence the general constructions behind the paper.
For technical reasons we actually prefer to work dually with dual
quasi-Hopf algebras $(kG,\phi)$ where $kG$ is the group algebra,
although this is equivalent when $G$ is finite. We introduce in
particular the example $k_FG$ as the quasi-associative algebra
associated to the group algebra of $G$ by twisting. Its algebraic
properties are studied in Section~3. Then Section~4 presents the
octonions (and more trivially, the quaternions and the complex
numbers) as examples of this type. We show how many of their
properties may be understood in terms of the 2-cochain $F$.
Section~5 provides new quasi-associative algebras beyond the
octonions. Section~6 introduces a suitable quasi-Hopf algebra of
`automorphisms' or comeasurings associated to any quasi-algebra of
the type we consider. Finally, Section~7 develops some first steps
in quasi-linear algebra, meaning a theory of matrix representations
of quasi-algebras.

\section{Preliminaries: General constructions}

An introduction to quantum groups, including quasitensor categories
and quasi-Hopf algebras is in \cite{Ma:book}, the main notations of
which we use here. In fact, the natural setting for us is the dual
of Drinfeld's axioms\cite{Dri:qua}, namely the notion of a dual
quasi-Hopf algebra\cite{Ma:tan}\cite{Ma:book}.

Thus, a {\em dual quasibialgebra} is a $(H,\Delta,\eps,\phi)$ where
the coproduct $\Delta:H\to H\tens H$ and counit $\eps:H\to k$ form
a coalgebra (the axioms are those of a unital associative algebra
with arrows reversed) and are multiplicative with respect to a
`product' $H\tens H\to H$. This is required to be associative up
to`conjugation' by $\phi$ in the sense
 \eqn{qassoc}{ \sum a\o\cdot(b\o\cdot c\o)\phi (a\t, b\t, c\t)
 =\sum \phi (a\o, b\o, c\o)(a\t\cdot b\t)\cdot c\t}
for all $a,b,c\in H$. Here $\Delta h=\sum h\o\tens h\t$ is a
notation and $\phi$ is a unital 3-cocycle in the sense
\eqn{3-cocycle}{\sum \phi(b\o, c\o, d\o)
\phi(a\o, b\t c\t, d\t)
\phi(a\t, b\th, c\th) =\sum\phi(a\o,b\o,c\o d\o)
\phi(a\t b\t,  c\t,  d\t),}
for all $a,b,c,d\in H$, and $\phi(a,1,b)=\eps(a)\eps(b)$ for all
$a,b\in H$. We also require that $\phi$ is convolution-invertible
in the algebra of maps $H^{\tens 3}\to k$, i.e. that there exists
$\phi^{-1}:H^{\tens 3}\to k$ such that
\[ \sum \phi(a\o,b\o,c\o)\phi^{-1}(a\t,b\t,c\t)=\eps(a)\eps(b)\eps(c)
=\sum \phi(a\o,b\o,c\o)\phi^{-1}(a\t,b\t,c\t)\]
for all $a,b,c\in H$.

A dual quasibialgebra is a quasi-Hopf algebra if there is a linear
map $S:H\to H$ and linear functionals $\alpha,\beta:H\to k$ such
that
\eqn{qcoant1}{\sum (Sa\o) a\th \alpha(a\t)=1\alpha(a),\quad
\sum a\o Sa\th \beta(a\t)=1\beta(a),}
\eqn{qcoant2}{\sum \phi (a\o, Sa\th, a\fiv)\beta(a\t)
\alpha(a\fo)=\eps(a),\quad
\sum\phi ^{-1}(Sa\o, a\th, Sa\fiv)\alpha(a\t)\beta(a\fo)
=\eps(a)}
for all $a\in H$.

Finally, $H$ is called {\em dual quasitriangular} if there is a
convolution-invertible map $\CR:H\tens H\to k$ such that
\eqn{qcoqua1}{\CR(a\cdot b, c) =\sum \phi (c\o, a\o,
b\o)\CR(a\t, c\t)
\phi^{-1}(a\th,c\th,b\t)\CR(b\th, c\fo)\phi (a\fo, b\fo, c\fiv),}
\eqn{qcoqua2}{\CR(a,b\cdot c)=\sum \phi ^{-1}(b\o, c\o, a\o)
\CR(a\t, c\t)\phi (b\t, a\th, c\th)\CR(a\fo, b\th)\phi ^{-1}
(a\fiv, b\fo, c\fo),}
\eqn{qcoqua3}{\sum b\o\cdot a\o\CR(a\t, b\t)
=\sum \CR(a\o, b\o)a\t\cdot b\t}
for all $a,b,c\in H$.

We recall also that a {\em corepresentation} or comodule under a
coalgebra means vector space $V$ and a map $\beta:V\to V\tens H$
obeying $(\id\tens\Delta)\circ\beta=(\beta\tens\id)\circ\Delta$ and
$(\id\tens\eps)\circ\beta=\id$. This is the notion of an action
with arrows reversed. In the finite-dimensional case a coaction of
$H$ means an action of the associative algebra $H^*$.

A monoidal category is a category $\CC$ of objects $V,W,Z$,etc. a
functor $\tens:\CC\times\CC\to \CC$ and a natural transformation
$\Phi: ((\ \tens\ )\tens\ )\to (\ \tens(\ \tens\ ))$ between the
two functors $\CC\times\CC\times\CC\to \CC$, where $\Phi$ obeys Mac
Lane's `pentagon identity' for equality of the two obvious
isomorphisms\cite{Mac:cat}
\[  (((V\tens W)\tens Z)\tens U\isom V\tens(W\tens (Z\tens U)))\]
built from $\Phi$ for any four objects $V,W,Z,U$. A braided or
`quasitensor' category is a monoidal one which has, in addition, a
natural transformation $\Psi:\tens\to
\tens^{\rm op}$ obeying two `hexagon' coherence conditions, see
\cite{JoyStr:bra}.

The comodules $\CM^H$ over a quasi-Hopf algebra form such a
category with
\[ \Phi_{V,W,Z}((v\tens w)\tens z)=\sum v\bo\tens( w\bo \tens z\bo)
\phi(v\bt,w\bt,z\bt)\]
for $v\in V$, $w\in W$, $z\in Z$. Here $\beta(v)=\sum v\bo\tens
v\bt$ is a notation and the tensor product is two comodules is a
comodule by composition with the product of $H$. In the
quasitriangular case the category is braided\cite{Dri:qua}, with
$\Psi_{V,W}(v\tens w)=\sum w\bo\tens v\bo \CR(v\bt,w\bt)$. There is
also a conjugate or dual coaction on $V^*$ made possible by the
antipode $S$. The converse is also true, namely any (braided)
monoidal category with duals and with a multiplicative functor to
the category of vector spaces (and some finiteness properties)
comes from the comodules over a dual (quasitriangular) quasi-Hopf
algebra\cite{Ma:tan}.

If $H$ is a dual quasi-Hopf algebra then so is $H_F$ with the new
product, $\Phi,\CR,\alpha,\beta$ given by
\ceqn{twist}{ a\cdot_F b=\sum F^{-1}(a\o,b\o)a\t b\t F(a\th,b\th)\\
\phi_F(a,b,c)=\sum F^{-1}(b\o,c\o) F^{-1}(a\o,b\t c\t)
 \phi(a\t,b\th,c\th)F(a\th b\fo,c\fo) F(a\fo,b\fiv)
 \\
\alpha_F(a)= \sum F(Sa\o,a\th)\alpha(a\t)
,\quad \beta_F(a)=\sum F^{-1}(a\o,Sa\th)\beta(a\t)\\
\CR_F(a,b)=\sum F^{-1}(b\o,a\o)\CR(a\t,b\t) F(a\th, b\th)}
for all $a,b,c\in H$. Here $F$ is any  convolution-invertible map
$F:H\tens H\to k$ obeying $F(a,1)=F(1,a)=\eps(a)$ for all $a\in H$
(a 2-cochain). This is the dual version of the twisting operation
or `gauge equivalence' of Drinfeld, so called because it does not
change the category of comodules up to monoidal equivalence.

\begin{defin} Let $H$ be a dual quasi-Hopf algebra. An $H$-comodule
quasialgebra $A$ is an algebra in the category of $H$-comodules.
This means an $H$-comodule, a product map $\cdot$ associative in
the category and equivariant under the action of $H$. Explicitly,
\[  (a\cdot b)\cdot c=\sum a\bo\cdot(b\bo\cdot c\bo) \phi(a\bt,b\bt,c\bt),
\quad  \beta(a\cdot b)= \beta(a)\beta(b) ,\quad\forall a,b\in H \]
where the last expression uses the tensor product algebra in
$A\tens H$.
\end{defin}

\begin{propos} If $A$ is an $H$-comodule quasialgebra and $F:H\tens H\to k$
a 2-cochain, then $A_F$ with the new product
\[ a\cdot_F b= \sum a\bo b\bo F(a\bt,b\bt)\]
and unchanged unit, is an $H_F$-comodule quasialgebra.
\end{propos}
\proof This is elementary and follows from the equivalence of the comodule
categories under twisting. See also \cite{GurMa:bra} in the module
version.
\endproof

There is a parallel theory with all arrows reversed. Thus, $H$ can
be a quasi-Hopf algebra, with associative product and $\Delta$
coassociative up to conjugation by an invertible 3-cocycle $\phi\in
H^{\tens 3}$\cite{Dri:qua}. In this case the modules of $H$ form a
monoidal category and, in the quasitriangular case a braided one.
In this case we work with $H$-module quasialgebras and their
twistings by $F\in H\tens H$.

When the theory is developed in this comodule form, it is an easy
matter to specialize to the following class of examples: let $H=kG$
the group algebra of a group. This has coproduct etc.
\[ \Delta x=x\tens x,\quad \eps x=1,\quad Sx=x^{-1},\quad \forall x\in G\]
forming a Hopf algebra. However, for any point-wise invertible
group cocycle $\phi:G\times G\times G\to k$ in the sense
\eqn{3cocy}{ \phi(y,z,w)\phi(x,yz,w)\phi(x,y,z)=\phi(x,y,zw)\phi(xy,z,w),
\quad \phi(x,e,y)=1}
extended linearly to $kG^{\tens 3}$, we can regard $(kG,\phi)$ as a
dual quasi-Hopf algebra. Group inversion provides an antipode with
$\alpha=\eps, \beta(x)=1/\phi(x,x^{-1},x)$. Finally, a dual
quasitriangular structure is possible only when $G$ is Abelian and
corresponds to invertible $\CR:G\times G\to k$ such that
\eqn{bichar}{ \CR(xy,z)=\CR(x,z)\CR(y,z){\phi(z,x,y)\phi(x,y,z)\over
\phi(x,z,y)},\quad
\CR(x,yz)=\CR(x,z)\CR(x,y){\phi(y,x,z)\over\phi(y,z,x)\phi(x,y,z)}}
for all $x,y,z\i G$.

A special case is when $\phi$ is a coboundary
\eqn{cobound}{ \phi(x,y,z)={F(x,y)F(xy,z)\over F(y,z)F(x,yz)},
\quad \CR(x,y)=\CR_0(x,y){F(x,y)\over F(y,x)},\quad \beta(x)
={F(x^{-1},x)\over F(x,x^{-1})}}
for any invertible $F$ obeying $F(x,e)=1=F(e,x)$ for all $x\in G$
and any invertible bicharacter $\CR_0$. This is the twisting of the
group algebra $(kG,\CR_0)$ regarded as a dual quasitriangular Hopf
algebra\cite{Ma:book} with trivial initial $\phi_0$.

Next, a coaction of $kG$ means precisely a $G$-grading, where
$\beta(v)=v\tens |v|$ on homogeneous elements of degree $|v|$.
Hence the notion of an $H$-comodule quasialgebra in this case
becomes:

\begin{defin} A $G$-{\em graded quasialgebra} is an $G$-graded
vector space $A$, a product map $A\tens A\to A$ preserving the
total degree and associative in the sense
\[ (a\cdot b)\cdot c =a\cdot (b\cdot c)\phi(|a|,|b|,|c|),
\quad\forall a,b,c\in A\]
(of homogeneous degree), for a 3-cocycle $\phi$. A $G$-graded
quasi-algebra is called {\em coboundary} if $\phi$ is a coboundary
as in (\ref{cobound}).
\end{defin}

This is the setting which we will use, with the above as the
underlying explanation of the constructions. If $G$ is finite we
can equally regard its function algebra $k(G)$ with $\phi\in
k(G)^{\tens 3}$ as a quasi-Hopf algebra $k_\phi(G)$ as in
\cite{Ma:qdou}, and then view a $G$-graded quasialgebra as
equivalently a $k_\phi(G)$-module quasialgebra. Here the action of
$h\in k_\phi(G)$ is $h.v=v h(|v|)$ on homogeneous elements.

\begin{corol} $k_FG$ defined as $kG$ with a modified product
\[ x\cdot_F y=xy F(x,y),\quad \forall x,y\in G\]
is a coboudary $G$-graded quasialgebra. The degree of $x\in G$ is
 $x$, and $F$ is any 2-cochain on $G$.
\end{corol}
\proof This is a special case of the twisting proposition.
Here $kG$ coacts on itself by $\beta=\Delta$, i.e. the degree of
$x\in G$ is $x$. We now twist $kG$ to the dual quasi-Hopf algebra
$(kG,\phi=\del F)$. In the process, we also twist $kG$ as a
comodule algebra to $k_FG$ as a comodule quasialgebra under this
dual quasi-Hopf algebra. \endproof

We note an elementary properties of $k_FG$.

\begin{propos} With trivial initial bicharacter $\CR_0$, the category
of $G$-graded vector spaces is symmetrically braided (in the sense
$\Psi^2=\id$) and $k_FG$ is braided-commutative in the sense
\[ a\cdot b=\cdot\circ\Psi(a\tens b),\quad\forall a,b\in k_FG.\]
\end{propos}
\proof This is immediate from (\ref{cobound}) and the definition of $\Psi$
from $\CR$. The latter is clearly $\Psi(a\tens b)=b\tens a
{F(|a|,|b|)\over F(|b|,|a|)}$ for elements of homogeneous  degree
$|a|,|b|$, which is braided but trivially braided in the sense
$\Psi^2=\id$ (i.e. the category of $G$-graded spaces in this case
is symmetric monoidal rather than strictly braided.) \endproof

\section{More about the quasi-algebras $k_FG$}

In this section, we will study further properties of the $G$-graded
quasialgebras $k_FG$ beyond the general ones arising from their
categorical structure in the preceding section. We assume that $G$
is Abelian, that $F$ is a 2-cochain and $\phi$ a 3-cocycle.

First of all, we note that $k_FG$ has a natural symmetric bilinear
form whereby the basis of group elements is orthonormal. In
general, the associated quadratic function on $k_FG$ will not be
multiplicative (a quadratic character).

\begin{propos} If $k_FG$ admits a quadratic character then $F^2$ is a
coboundary and $\phi^2=1$. If the Euclidean norm quadratic function
defined by $q(x)=1$ for all $x\in G$ is multiplicative (making
$k_FG$ a composition algebra) then $F^2=1$.
\end{propos}
\proof (For all
discussions of quadratic forms we suppose that $k$ has
characteristic not $2$). Given a quadratic character $q:k_FG\to k$,
we have $q(x\cdot y)=q(F(x,y)xy)=F^{2}(x,y)q(xy)=q(x)q(y)$ on
$x,y\in G$, i.e. $F^2(x,y)=q(x)q(y)/q(xy)$ is a coboundary in the
group cohomology. In general, if $F^2=\del q$ we still need to
specify a bilinear form with diagonal $q$, so the converse is not
automatic. If we take the canonical quadratic function associated
to $G$ as an orthonormal basis, we will have
$q(x\cdot_Fy)=F^2(x,y)=F^2(x,y)q(x)q(y)$ for all $x,y\in G$, so if
this is multiplicative then $F^2=1$.  \endproof

This will be the case for some of the Cayley algebras in the next
section, as well as for many other examples, and is the reason that
$F,\phi$ typically have values $\pm1$ in these cases.

Also, we already know from the construction in Corollary~2.4 that
\[(x\cdot y)\cdot z=x\cdot (y \cdot z)\phi(x,y,z),\quad
\phi(x,y,z)={F(x,y)F(xy,z)\over F(x,yz)F(y,z)}\]
and hence that $k_FG$ is associative {\em iff} $\phi=1$. Also,
recall that $F(x,e)=F(e,x)=1$ (where $e\in G$ is the group
identity) is part of the cochain definition, and \[
\phi(e,x,y)=\phi(x,e,y)=\phi(x,y,e)=1\] is part of the cocycle
definition (the middle one implies the other two), in particular it
holds for our coboundary $\phi$.

Likewise, we know from Proposition~2.5 that $k_FG$ is
braided-commutative with respect the braiding
\[ \Psi(x\tens
y)=\CR(x,y)y\tens x,\quad \CR={F(x,y)\over F(y,x)},\quad\forall
x,y\in G.\] Hence it is commutative in the usual sense {\em iff}
$F$ is symmetric. This is also clear from the form of the product
in $k_FG$ since $G$ itself is Abelian. More interesting for us,

\begin{defin} We say that $k_FG$ is {\em altercommutative} if
$\Psi$ is given by $\CR$ of the form
\[ \CR(x,
y)=\cases{1&if\ $x=e$\ or\ $y=e$\ or\ $x=y$\cr -1& otherwise\cr}\]
for all $x,y\in G$.
\end{defin}

Note that $\CR(x,x)=\CR(x,e)=\CR(e,x)=1$ for any $k_FG$, so the
content here is the value $-1$ in the remaining `otherwise' case.
Also note that an altercommutative $k_FG$ can never be commutative
unless $G=\Z_2$.  The condition is somewhat similar to the notion
of a `supercommutative' algebra. One also has (more familiar) cases
for the breakdown of associativity, such as the notion of an
alternative algebra. We have,

\begin{propos} \kfg\ is an alternative algebra if and only if
\[ \phi^{-1}(y,x,z) + \CR(x,y)\phi^{-1}(x,y,z)=1+\CR(x,y)\]
\[ \phi(x ,y  ,z   )+{\CR(z,y)}\phi(x ,z   ,y  )=1+{\CR(z,y)}\]
for all $x,y,z\in G$. In this case,
\[  \phi(x,x,y)=\phi(x,y,y)=\phi(x,y,x)=1 \] for all
$x,y\in G$.
\end{propos}

\proof It is enough to consider the conditions of an alternative
algebra on our basis elements, $x,y,z\in G$, i.e.
\cmath{(x\cdot
y)\cdot z-x\cdot (y\cdot z)+(y\cdot x)\cdot z-y\cdot (x\cdot z)=0\\
         (x\cdot y)\cdot z-x\cdot(y\cdot z)+(x\cdot z)\cdot y-x\cdot
(z\cdot y)=0.} This translates at once into the two equations
\[ F(x,y)
F(xy,z)+F(y,x) F(yx,z)= F(y,z) F(x,yz)+ F(x,z) F(y,xz)\] \[ F(x,z)
F(xz,y)+ F(x,y) F(xy,z)= F(z,y) F(x,zy)+ F(y,z) F(x,yz)\] for all
$x,y,z\in G$ when we put the product of $k_FG$ in terms of the
associative product in $G$. Dividing through then gives the
equations in terms of $\phi,\CR$ as stated.

Also, setting $x=y$ in the first equation gives us (for
characteristic of $k$ not 2) $\phi(x,x,z)=1$. Setting $y=z$ in the
second equation likewise gives us $\phi(x,y,y)=1$. Given these,
setting $x=z$ in either gives $\phi(x,y,x)=1$. Actually, it is
known that the condition of being an alternating algebra is
equivalent to
\[(
a\cdot a)\cdot b=a\cdot (a\cdot b),\quad (a\cdot b)\cdot
b=a\cdot(b\cdot b)\] for all $a,b$ in the algebra, which more
immediately implies $\phi(x,x,y)=\phi(x,y,y)=1$ on basis elements.
(Given these, the same two equations applied to $a=x+y, b=z$ in the
first case and $a=x, b=y+z$ in the second case provide the full
equations for an alternating algebra on basis elements $x,y,z$, and
hence imply that $(a\cdot b)\cdot a=a\cdot(b\cdot a)$ holds as
well, as usual.) \endproof

Next we consider involutions. Since we have a special basis of
$k_FG$ it is natural to consider involutions diagonal in this
basis.

\begin{lemma} \kfg\ admits an involution which is diagonal in the basis
$G$ {\em iff} $\CR=\del s$ (a group coboundary) for some 1-cochain
$s:G\to k^*$ with $s^2=1$. In this case, one has
$\CR(x,y)=\CR(y,x)$ and $\phi(x,y,z)=\phi(z,y,x)^{-1}$ for all
$x,y,z\in G$.
\end{lemma} \proof
Consider the endomorphism $\sigma$ of the vector space \kfg\ of the
form \[
\sigma(x)=s(x)x,\quad \forall x\in G,\] extended linearly,  for some
function $s:G\to k$. For an involution we need $\sigma^2=\id$,
$s(1)=1$ and $\sigma(a \cdot b)=\sigma(b)\sigma(a),\ \forall a,b$
in our algebra. It is enough to consider these on the basis
elements. Then clearly the first two correspond to

i) $s(e)=1$ and $s^2(x)=1$ for all $x\in G$.

For the second condition,
\cmath{\sigma(x\cdot y)=s(xy) x\cdot y=s(xy)F(x,y) xy}
while,
\[\sigma(y)\cdot\sigma(x)=s(y)y \cdot
s(x)x=s(x)s(y) F(y,x)xy.\] Equality for all $x,y$ corresponds to

ii) ${s(x)s(y)\over s(xy)}= {F(x,y)\over F(y,x)}$ for all $x,y \in
G$.

We interpret this as stated, where the right hand side is $\CR$
corresponding to the braiding $\Psi$. In view of i), it implies
that $\CR^2(x,y)=1$ and hence that $\CR$ is symmetric. It also
implies that $\del\CR=1$ in the group cohomology, which is the
condition on $\phi$ stated since $\CR(x,y)=F(x,y)/F(y,x)$ for all
$x,y\in G$. \endproof

In fact, we will be particularly interested in involutions with the
property,
\[ a+\sigma(a), a\cdot \sigma(a)\in k 1\]
(a multiple of the identity) for all $a$ in the algebra. Let us
call this a {\em strong involution}.

\begin{propos}   $k_FG$ admits a diagonal strong involution $\sigma$ if
and only if

i) $G\simeq (\Z_2)^n$ for some $n$,

ii)  $\sigma(e)=e$, $\sigma(x)=-x$ for all $x\ne e$,

iii) $k_FG$ is altercommutative in the sense of Definition~3.2.

\end{propos} \proof Consider an endomorphism of the diagonal form
$\sigma(x)=s(x)x$. A  general element $a=\sum_{x\in G} \alpha_x x$,
we have of course
\[ \sigma(a)= \sum_{x\in G} \alpha_x s(x)x,\] hence
\[ a+\sigma(a) =\sum_{x\in G}\alpha_x(1+s(x))x\in k1\]
for all coefficients $\alpha_x$ if and only if $s(x)=-1$ for all
$x\ne e$. Here the group identity $e\in G$ is the basis element
$1\in k_FG$. Since we also need $\sigma(1)=1$ for an involution,
these two conditions hold iff

ii') $s(e)=1$ and $s(x)=-1$ for all $x\ne e$.

\noindent Next, consider a basis element $x\in G$, then \[
x\cdot\phi(x)=s(x)x\cdot x=s(x)F(x,x)x^2\in k1\] implies that
$x^2=e$, i.e. every element of $G$ has order 2. Hence $G$ (which is
a finite Abelian group) is isomorphic to $(\Z_2)^n$ for some $n$.
Next, consider $x+y$ for basis elements $x\ne y$.  Then
\[(x+y)\cdot\phi(x+y)=(x+y)\cdot(s(x)x+s(y)y)\qquad\]
\[\qquad=(s(x)F(x,x)+s(y)F(y,y))e +(s(x)F(x,y)+s(y)F(y,x))xy\in k1\]
tells us that

iii') $s(x)F(y,x)+s(y)F(x,y)=0$ for all $x\not= y$.

When $x=e$ or $y=e$ this is empty. Otherwise, given  ii') it is
equivalent to $\CR(x,y)=-1$ for all $x\ne y$, $x\ne e$ and $y\ne
e$, which is the altercommutativity condition.

Conversely, given these facts, a general element $a=\sum_{x\in
G}\alpha_x x$ obeys
\[a\cdot\sigma(a)=\sum_{x,y\in G}\alpha_x\alpha_y
s(y)x\cdot y=\sum_{x,y\in G}\alpha_x\alpha_y F(x,y)s(y)xy\qquad\]
\[=\sum_{z\in G}z\left(\sum_{x,y\in G}\delta_{xy,z}\alpha_x\alpha_y
s(y)F(x,y)\right)\in k1\] since the terms with $z\ne e$ have
contribution only when $x\ne y$, and in this case the terms cancel
pairwise due to condition iii').  Hence if these conditions hold,
we have the stated properties for $\sigma$.  They clearly imply the
ones in the preceding lemma as well. \endproof

On the other hand, given $k_FG$ we have a natural function
$s(x)=F(x,x)$ and consider now the particular endomorphism
corresponding to this.

\begin{propos} If $\sigma(x)=F(x,x)x$ for all $x\in G$ defines a strong
involution and $|G|\ne 2$ in $k$, then the algebra \kfg\ is simple.
\end{propos}
\proof Let $I$ be an ideal of \kfg\ different from \kfg,
and $a=\sum_{x\in G}\alpha_x x $ an element of $I$. Then for each
$y\ne e\in G$,
\[ a\cdot y =\sum_{x\in G}\alpha_x F(x,y)xy \in I,\quad y
\cdot a =\sum_{x\in G}\alpha_x F(y,x)xy \in I.\]
Adding these together and using Proposition~3.5 in the form of the
conditions ii') and iii'), we see that
\[a\cdot y+y\cdot a=\sum_{x\ne y}\alpha_x(F(y,x)+F(x,y))xy
+\alpha_y2F(y,y)e=2\alpha_ey-2\alpha_ye\in I.\] Then $(a\cdot
y+y\cdot a)\cdot y=-2(\alpha_ee+\alpha_y y)\in I$ as well. Since
this is for all $y\ne e$, we have $\sum_{y\ne
e}\alpha_yy+(|G|-1)\alpha_e e\in I$ and hence $\alpha_ee\in I$
provided $|G|-2\ne 0$ in $k$. In this case $\alpha_e=0$ and
$-\alpha_y e+\alpha_e y=-\alpha_ye\in I$ tells us that $\alpha_y=0$
for all $y$, i.e. $a=0$. Hence $I=\{0\}$.  \endproof

We also see from Proposition~3.5 that $\sigma$ a strong involution
restricts us to $G\isom(\Z_2)^n$. In this context we have a partial
converse to Proposition~3.1.

\begin{propos} If $G\isom(\Z_2)^n$ then the Euclidean norm quadratic
function defined by $q(x)=1$ for all $x\in G$ makes $k_FG$ a
composition algebra if and only if

i) $F^2(x,y)=1$ for all $x,y\in G$

ii) $F(x,xz)F(y,yz)+F(x,yz)F(y,xz)=0$, for all $x,y,z\in G$ with
$x\ne y$.

In this case the conditions in the Proposition~3.6 hold, i.e.
$\sigma(x)=F(x,x)x$ for all $x\in G$ is a strong involution and
$k_FG$ is simple if $|G|\ne 2$. Moreover, $k_FG$ is alternative.
\end{propos}

\proof Suppose that $q$ is multiplicative.  $q(x\cdot y)=q(x)q(y)$ on
all basis elements $x,y\in G$ is $F^2=1$, as we know already from
Proposition~3.1. The next case $q((x+y)\cdot z)=q(x+y)q(z)$ on
basis elements $x,y,z\in G$ with $x\ne y$ does not yield any new
condition (both sides are 2). Now consider the elements $x,y,z,w\in
G$ with $x\ne y$, $z\ne w$ but $xz=yw$. Because every element of
$G$ is of order 2, this also means $xw=yz$. Because $G$ is a group,
$xz\ne xw$, however. Hence
\[
q((x+y)\cdot(z+w))=q(xz(F(x,z)+F(y,t))+xw(F(x,w)+F(y,z)))\quad\]
\[\qquad=(F(x,z)+F(y,w))^2
+(F(x,w)+F(y,z))^2=4+2(F(x,z)F(y,w)+F(x,w)F(y,z))\] while
$q(x+y)q(z+w)=4$. This is the second condition stated after writing
$w=xyz$ and renaming $xz$ to $z$.

Conversely, assuming these conditions and given $a=\sum_x\alpha_x
x$ and $b=\sum_y\beta_y y$, we have
\[ a\cdot b=\sum_{z\in
G}z\left(\sum_{x\in G}\alpha_x\beta_{xz}F(x,xz)\right)\] since
every element of $G$ has order 2. Hence
\[
q(ab)=\sum_z\left(\sum_x\alpha_x\beta_{xz}F(x,xz)\right)^2=
\sum_z\sum_x\sum_y
\alpha_x\alpha_y\beta_{xz}\beta_{yz}F(x,xz)F(y,yz).\]
In this sum the diagonal part where $x=y$ contributes
\[ \sum_z \sum_x
\alpha_x^2\beta_{xz}^2=q(a)q(b)\]
since $F^2=1$, while the remaining contribution from $x\ne y$ is
\[ \sum_z\sum_{x\ne y} \alpha_x\alpha_y\beta_{xz}\beta_{yz}F(x,xz)
F(y,yz),\quad\qquad (*)\] By condition ii) this is equal to
\[ -\sum_{x\ne
y}\sum_z\alpha_x\alpha_y\beta_{xz}\beta_{yz}F(x,yz)F(y,xz)
=-\sum_{x\ne y}\sum_{w}
\alpha_x\alpha_y\beta_{yw}\beta_{xw}F(x,xw)F(y,yw)\] where we
change the order of summation and change variables to $w=xyz$. But
this has the same form as our original expression for (*) but with
a minus sign, hence this is zero.

Next, we observe that the condition ii) can be broken down
equivalently as the conditions

ii.a) $F(x,xy)+F(x,y)=0$ for all $x,y\in G$ with $x\ne e$.

ii.b) $F(x,yz)F(y,xz)+F(x,z)F(y,z)=0$ for all $x,y,z\in G$ with
$x\ne e$, $y\ne e$ and $x\ne y$.

The first of these is ii) in either of the cases $x=e,y\ne e$ or
$x\ne e,y=e$ (followed by a relabeling), while the second is the
remaining case $x\ne e,y\ne e,x\ne y$ after making use of ii.a) to
substitute $F(x,xz)$ and $F(y,yz)$.

In this case, ii.a) implies the condition ii') of Proposition~3.6.
On the other hand, $z=e$ in the original form of the present
condition ii) gives us \[ F(x,x)F(y,y)+F(x,y)F(y,x)=0,\quad \forall
x\ne y\] which, given $F^2=1$, implies the condition iii') of
Proposition~3.6. Hence $\phi$ is strongly involutive.

Finally, also in this case, the equations of an alternative algebra
in terms of $F$ (see the proof of Proposition~3.3) reduce to the
following. If $x=e$ or $y=e$, the first equation is trivial.
Otherwise, the case $x\ne e, y\ne e, x=y$ reduces to
$\phi(x,x,z)=1$, which in our present case where $G\isom(\Z_2)^n$
reduces to
\[ F(x,x)=F(x,xz)F(x,z).\]
This holds because the left hand side is -1 by the condition ii')
of Proposition~3.6 and the right hand side is $-F(x,z)^2=-1$ by
ii.a) and $F^2=1$. The remaining case is $x\ne e,y\ne e$ and $x\ne
y$. In this case the altercommutativity property in Proposition~3.6
reduces the first equation for an alternative algebra to
\[ F(y,z)F(x,yz)+F(x,z)F(y,xz)=0.\]
Since $F^2=1$, this is equivalent to ii.b)

On the other hand, under the assumptions of Proposition~3.6, the
conditions ii.a)-ii.b) are equivalent to

ii.c) $F(x,y)+F(xy,y)=0$ for all $x,y\in G$ with $y\ne e$

ii.d) $F(x,y)F(x,z)+F(xy,z)F(xz,y)=0$ for all $x,y,z\in G$ with
$y\ne e$, $z\ne e$ and $y\ne z$

(or one can obtain them directly from our original condition ii)).
We use these versions in a similar analysis for the content of the
second of the conditions  in Proposition~3.3 for an alternative
algebra.
\endproof

We have written the proof of the last part of the proposition in a
reversible way. Hence we also conclude,

\begin{corol} If $\sigma(x)=F(x,x)x$ for all $x\in G$ is a strong
involution, and $F^2=1$, then the following are equivalent,

i) $k_FG$ is an alternative algebra,

ii) $k_FG$ is a composition algebra.
\end{corol}

The conditions in Proposition~3.7 and the corollary are evidently
highly restrictive, because if $k$ has characteristic different
from 2 it is known that we have only the following composition
algebras with the Euclidean norm: $k$, the algebra with basis
$1,v,v^2=-1$, the algebra over $k$ with the product of quaternions
and the algebra over $k$ with the product of octonions. If $k$ is
algebraically closed these algebras are isomorphic to $k$, $k\oplus
k$, $M_2(k)$ ($2\times 2$ matrices) and Zorn's algebra of vectorial
matrices \cite{ZSSS:rin}.  Equivalently, one knows that these are
the only simple alternative algebras. On the other hand,  the
diagonal strong involution conditions in Proposition~3.6 are
definitely weaker and hold for the entire family of Cayley
algebras, as we will see in the next section.

Finally, for completeness, we include a slight generalisation of
Proposition~3.7.

\begin{propos} Let $k_FG$ be an algebra that admits a strong diagonal
involution $\sigma (x)= s(x)x$. Then the non degenerate form
$n(x)=x\cdot \sigma (x)$ makes $k_FG$ a composition algebra if and
only if

i) $s(xy)F(x,y)^2 F(xy,xy)=s(x)s(y)F(x,x)F(y,y)$, for all $x,y\in
G$.

ii) $F(x,xz)F(y,yz)F(z,z)s(z)+F(x,yz)F(y,xz)F(xyz,xyz)s(xyz)=0$,
for all $x,y\in G$ with $x\not=y$.
\end{propos}
\proof Let $A$ be a quasialgebra $k_FG$ that admits a strong diagonal
involution $\sigma (x)=s(x)x$ and let us consider the form
$n(x)=x\sigma(x)$. For all $x,y\in G$ if the form $n(x)$ admits
composition we have  $n(x\cdot y)=n(x)n(y)$ and then $n(x\cdot
y)=F^2(x,y)F(xy,xy)s(xy)=n(x)n(y)=s(x)s(y)F(x,x)F(y,y)$. On the
other hand for two  elements of the algebra $k_FG$, $a=\sum_{x\in
G}\alpha_x x$ and $b=\sum_{y\in G}\beta_y y$ we have
$n(a)=\sum_{x\in G}\alpha^2_x s(x)F(x,x)$ and $n(b)=\sum_{y\in
G}\beta^2_y s(y)F(y,y)$. But we know that $a\cdot b=\sum_{z\in G}
z(\sum_{x\in G}\alpha_x
\beta_{xz}F(x,xz))$ so $n(a\cdot b)=\sum_{z\in G}[\sum_{x\in G}\alpha_x
\beta_{xz}F(x,xz)]^2 s(z) F(z,z)=\sum_{x,y,z\in
G}\alpha_x\beta_{xz}F(x,xz)\alpha_y\beta_{yz}F(y,yz) s(z) F(z,z)$,
and like in the proof of Proposition~3.7 the result follows after
comparing the last expression with
\[ n(a)n(b)=\left(\sum_{x\in G}\alpha^2_x s(x)F(x,x)\right)
\left(\sum_{y\in G}\beta^2_y
s(y)F(y,y)\right).\]
\endproof

\section {Cayley algebras}

In this section, we show that the `complex number' algebra, the
quaternion algebra, the octonion algebra and the higher Cayley
algebras are all $G$-graded quasialgebras the form $k_FG$ for
suitable $G$ and $F$, which we construct. We recall that these
algebras can be constructed inductively by the Cayley-Dickson
process; we show that this process is compatible with our
quasialgebra approach to nonassociative algebras.

Let  $A$ be finite-dimensional (not necessarily associative)
algebra with identity element 1 and a strong involution $\sigma$,
i.e. an involution such that $a+\sigma(a),a\cdot\sigma(a)\in k1$
for all $a\in A$. We have studied this condition in the context of
our quasialgebras $k_FG$ in Proposition~3.5. The Cayley Dickson
process says that we can obtain a new algebra $\bar A=A\oplus vA$
of twice the dimension (i.e. elements are denoted $a,va$ for $a\in
A$) and multiplication defined by
\[
(a+vb)\cdot (c+vd)=(a\cdot c+\alpha d\cdot
\sigma(b))+v(\sigma(a)\cdot d+ c\cdot b) \]
 and with a new
strong involution $\bar\sigma$
\[ \bar\sigma(a+vb)=\sigma(a)-vb.\]
The symbol $v$ here is a notation device to label the second copy
of $A$ in $\bar A$. However, $v\cdot v=\alpha 1$ according to the
stated product, so one should think of the construction as a
generalisation of the idea of complexification when $\alpha=-1$. If
$A=k$ and $\alpha=-1$ then $\bar A=k[v]$ modulo the relation
$v^2=-1$ will be called the `complex number algebra' over a general
field $k$. As in the preceding section, we suppose $k$ has
characteristic not 2.

We start with the cochain version of the Cayley-Dickson
construction, motivated by the formulae above.

\begin{propos} Let $G$ be a finite Abelian group $F$ a cochain on it
(so $k_FG$ is a $G$-graded quasialgebra). For any $s:G\to k^*$ with
$s(e)=1$ we define $\bar G=G\times \Z_2$ and on it the cochain
$\bar F$ and function $\bar s$,
\cmath{\bar F(x,y)=F(x,y),\quad \bar
F(x,vy)=s(x)F(x,y),\quad \bar F(vx,y)=F(y,x)\\ \bar F(vx,vy)=\alpha
s(x)F(y,x),\quad \bar s(x)=s(x),\quad \bar s(vx)=-1} for all
$x,y\in G$. Here $x\equiv(x,e)$ and $vx\equiv(x,\nu)$ denote
elements of $\bar G$, where $\Z_2=\{e,\nu\}$ with product
$\nu^2=e$.

If $\sigma(x)=s(x)x$ is a strong involution, then $k_{\bar F}\bar
G$ is the Cayley-Dickson process applied to $k_FG$.
\end{propos}
\proof The only features to be
checked for a cochain are that $\bar F$ should be pointwise
invertible (which is clear from invertibility of $s,F$) and $\bar
F(e,vx)=s(e)F(e,x)=1$ and $\bar F(vx,e)=F(e,x)=1$.  Hence we have a
new quasi-algebra $k_{\bar F}\bar G$.

This reproduces the product $\cdot$ of the Cayley-Dickson process
with respect to $\sigma(x)=s(x)x$, since that is
\cmath{ v x     \cdot y
=v(y\cdot  x )=vF(y,x)yx\\ v x\cdot vy      =\alpha y s(x)\cdot x
=\alpha s(x)
F(y,x)yx \\ x \cdot v y      =vs(x)x\cdot y =s(x) F(x,y)x v y} in
terms of the product of $\bar G$ on the right.  Moreover, $\bar s$
on $\bar G$ clearly reproduces the $\bar\sigma$ in the
Cayley-Dickson precedure as well. \endproof

This provides a cochain approach to the Cayley-Dickson process. The
last proposition also makes evident that all composition algebras
with identity element, over a field of characteristic different
from two, are in fact quasialgebras $k_FG$. We know that if we
start by a field $k$ with characteristic different from 2 we can
construct the sequence of algebras $k(\nu)$ (the Cayley-Dickson
extension of $k$ with $\alpha=\nu$), $(k(\nu),\beta)$ (the
Cayley-Dickson extension of $k(\nu)$ with $\alpha=\beta)$ and
$((k(\nu),\beta),\gamma)$ (the Cayley-Dickson extension of
$(k(\nu),\beta)$ with $\alpha=\gamma$), in all cases admitting a
strong diagonal involution  such that  the associated form admits
composition. If $A$ is a composition algebra with identity, over a
field of characteristic different from two, it is isomorphic to one
of these 4 classes of algebras. This proves too that the last
proposition of the section 3 is a complete characterization of all
composition algebras $k_FG$ with identity, in terms of its cochain
$F$. We also know from Proposition~3.5 that all these algebras with
strong diagonal involution are altercommutative.

As in Section~3, we are particularly interested in the canonical
involution defined by $s(x)=F(x,x)$ and in the Cayley Dickson
extension with $\alpha=-1$.
 In that case $\bar F$ is determined from $F$
alone. Let us call this choice the {\em standard Cayley-Dickson}
process in our cochain approach.

\begin{propos}

If $F(x,x)=-1$ for all $x\in G$ and $x\ne e$  then the same holds
for $\bar F$ under the standard Cayley-Dickson process. Moreover,
in this case,

(i) $\bar s$ has the standard form on $\bar G$.

(ii) If $F^2=1$ then $\bar F^2=1$ as well.

(iii) If $k_FG$ is altercommutative then so is $k_{\bar F}\bar G$.

\end{propos} \proof We have
\[
\bar\sigma(x)=\bar F(x,x)x=F(x,x)x=\sigma(x),\quad \bar\sigma(vx)=\bar
F(vx,vx)vx=-F(x,x)^2vx=-vx\] since $F^2(x,x)=1$. The first two
parts are then immediate. The third part is
\[\bar{\CR}(x,vy)={\bar F(x,vy)\over \bar F(vy,x)}=F(x,x),\quad
\bar{\CR}(vx,y)=F(y,y)^{-1}\]
\[ \bar{\CR}(vx,vy)={\bar F(vx,vy)\over F(vy,vx)}={F(x,x)\over
F(y,y)}\CR(x,y)\] has the altercommutative form in Definition~3.2
when $\CR$ does and when $F(x,x)=-1$ for $x\ne e$.
\endproof

In particular, the standard complex, quaternion, octonion etc
algebras are all of this form given by iterating the standard
Cayley-Dickson process. To describe their cochains, we consider the
special case where $G=(\Z_2)^n$ and $F$ is of the form
\[ F(x,y)=(-1)^{f(x,y)}\]
for some $\Z_2$-valued function $f$ on $G\times G$ (which is a
natural supposition for the class with $F^2=1$).

\begin{corol} If $G=(\Z_2)^n$ and $F=(-1)^f$ then the standard
Cayley-Dickson process has
 $\bar G=(\Z_2)^{n+1}$ and $\bar F=(-1)^{\bar f}$.  We use a vector
notation $\vec x=(x_1,\cdots,x_n)\in (\Z_2)^n$ where
$x_i\in\{0,1\}$ (and the group $\Z_2$ is now written additively).
Then
\[ \bar f((\vec
x,x_{n+1}),(\vec y,y_{n+1}))=f(\vec x,\vec y)(1-x_{n+1}) +f(\vec
y,\vec x)x_{n+1}+y_{n+1}f(\vec x,\vec x)+x_{n+1}y_{n+1}.\]
\end{corol}
\proof
>From the above, we have clearly \[ \bar f(x,y)=f(x,y),\quad \bar
f(x,vy)=f(x,y)+f(x,x),\] \[\bar f(vx,y)=f(y,x),\quad \bar
f(vx,vy)=1+f(x,x)+f(y,x).\] We then convert this to a vector
notation where each copy of $\Z_2$ in $G$ is the additive group of
$\Z_2$. We then make use of the {\em product} in $\Z_2$ to express
whether a term is included or not (thus $x_{n+1}y_{n+1}$
contributes 1 iff both $x_{n+1}=1$ and $y_{n+1}=1$, etc.).
\endproof

Iterating this now generates the $f$ for the quaternions, octonions
etc.:

\begin{propos}

i) The `complex number' algebra has this form with

 \[ G=\Z_2,\quad f(x,y)=xy,\quad x,y\in\Z_2\] where we identify $G$ as
the additive group $\Z_2$ but also make use of its product.

ii) The quaternion algebra is of this form with

\[ \bar G=\Z_2\times\Z_2,\quad \bar
f(\vec{x},\vec{y})=x_1y_1+(x_1+x_2)y_2\] where
$\vec{x}=(x_1,x_2)\in\bar G$ is a vector notation.

iii) The octonion algebra is of this form with

\[\bar{\bar G}=\Z_2\times\Z_2\times\Z_2,\quad \bar{\bar
f}(\vec{x},\vec{y})= \sum_{i\le
j}x_iy_j+y_1x_2x_3+x_1y_2x_3+x_2x_2y_3.\]

iv) The 16-onion algebra is of the form $\bar{\bar{\bar
G}}=\Z_2\times\Z_2\times\Z_2\times\Z_2$ and
\[\bar{\bar{\bar
f}}(\vec{x},\vec{y})=\sum_{i\le j}x_iy_j+ \sum_{i\ne j\ne k\ne
i}x_ix_jy_k+\sum_{{\rm distinct}\ i,j,k,l} x_ix_jy_ky_l+\sum_{i\ne
j\ne k\ne i}x_iy_jy_kx_4\]
\end{propos}

We are now able to apply our various criteria in the last section
for the structure of the algebras of the form $k_FG$, to this
construction and to all these algebras. Note in all these cases
(and for the who $2^n$-onion family generated in this way) $f$ has
a bilinear part defined by the bilinear form
\[ \pmatrix{1&1&\cdots &1\cr
0&1&\cdots&1\cr\vdots& & &\vdots\cr 0&\cdots&1&1\cr
0&\cdots&0&1}.\] For the complex number and quaternion algebras
this is the only part, which implies that $F$ is a bicharacter and
hence, in particular, its group coboundary $\phi=1$, i.e. these
algebras are associativity. The $f$ for the octonions has this
bilinear part,  which does not change associativity, plus a cubic
term which contributes to $\phi$. The $16$-onion has additional
cubic an quartic terms, etc. This makes the origin of the breakdown
of associativity for the higher members of the family particularly
clear.

In the remainder of this section we suppose that the quasialgebra
$k_FG$ is admits a diagonal involution $\sigma(x)=s(x) x$ in the
basis $G$ and we shall denote by $k_{\bar F}\bar G$ the generalized
Cayley-Dickson extension with respect to this
 (with general $\alpha \not=0$). It is easy to see that $\bar s$
 provides a
diagonal involution on it.

\begin{propos} The associativity
cocycle $\phi$ of $k_FG$ and the associativity cocycle $\bar\phi$
of $k_{\bar F}\bar G$ are related by
\cmath{\bar \phi     (x ,y  ,z   )=\phi(x ,y
,z   ),\quad \bar \phi     (v x     ,y  ,z   )={\CR(y,z)}\phi(x ,y
,z )\\ \bar \phi      (x ,v y  ,z )={\CR(y,z)}{\CR(xy,z)}\phi(x ,y
,z),\quad \bar \phi      (x ,y ,v z  )={\CR(x,y)}\phi(x ,y  ,z )\\
\bar \phi      (v x  ,v y ,z )={\CR(xy,z)}\phi(x   ,y,z ),\quad
\bar \phi   (v x   ,y  ,v z )={\CR(y,z)}{\CR(x,y)}\phi(x ,y ,z
)\\ \bar \phi (x ,v y ,v z      )=\CR(x,yz)\phi(x ,y  ,z ),\quad
\bar \phi (v x     ,v y       ,v z )={\CR(xy,z)}{\CR(x,y)}\phi(
x,y  ,z )} $\forall x,y,z\in G$.
\end{propos}
\proof We use the
definition of $\bar\phi$ as coboundary of $\bar F$, the form of
this in Proposition 4.1 and $F(x,y)={\CR(x,y)}F(y,x)$ from
Definition~3.2. We also use Lemma~3.4 which tells us that $s^2=1$,
$\CR(x,y)=s(x)s(y)/s(x,y)$ and $\phi(x,y,z)\phi(z,y,x)=1$.
\endproof

\begin{corol}
$k_{\bar F}\bar G$ is associative if and only if $k_FG$ is
associative and commutative.
\end{corol}
\proof We already know that $k_FG$ is commutative iff $F$ is symmetric,
which means iff ${\CR}=1$. So in this case if $\phi=1$ then
$\bar\phi=1$. So if $k_FG$ is associative and commutative, $k_{\bar
F}\bar G$ is associative. Conversely, if $\bar\phi=1$ then by
restriction, $\phi=1$ so $k_FG$ is associative. Moreover,
$1=\bar\phi(x,y,v)={\CR(x,y)}\phi(x,y,e)={\CR(x,y)}$ tells us that
$k_FG$ is commutative.  \endproof

\note{ If $\phi=1$ then Proposition~4.5 tells us the form of
$\bar\phi$ and we easily check that $k_{\bar F}\bar G$ is
alternative by Proposition~3.3. Thus,
\[\bar \phi     (v x     ,y  ,z   )+\bar{\CR}(z,y)\bar \phi(v x
,z   ,y
 )={\CR(y,z)}+1=\CR(z,y)+1=\bar{\CR}(z,y)+1$, etc by Lemma~3.4.
Converserly, suppose that $k_{\bar F}\bar G$ is alternative. We
know that $\phi(e,x ,y )=\phi(x ,e,z   )=\phi(x ,y  ,e)=1$ and as
$k_FG$ is alternative (because $k_{\bar F}\bar G (\alpha)$ is
alternative and the cocycle $\phi$ is a restriction of the cocycle
$\bar\phi$ etc.), so we have also $\phi(x ,x ,y )=\phi(x ,z   ,z
)=\phi(x ,y  ,x )=1$. So we have only to study the other cases.
But by Proposition~3.3, if $z\not=y,z\not=e,y\not=e$ we have
$\phi(x ,z   ,y  )=\phi(x ,y  ,z )$ and $\bar \phi (v x ,z   ,y
)=\bar \phi     (v x     ,y  ,z )$ which implies that $\phi(z   ,y
,x )=\phi(y  ,z   ,x )$.}

\section{New quasialgebras}

 It is clear that there are many examples of quasialgebras \kfg\
according to the group $G$ and the cochain $F$. In this section we
will consider examples where $F^2=1$; we will see that even in this
case we can obtain very different types of algebras.

By Proposition 3.7 we know that if \kfg\ is an $n$-dimensional
composition algebra  for the Euclidean norm, the cochain $F$ is
defined by an Hadamard matrix $H$ (that is, a matrix such that
$H^tH=nI$ where $t$ is transpose and $I$ is the identity matrix).
The matrix entries are $H_{x,y}=F(x,y)$ for $x,y\in G$. Motivated
by this, we begin by considering more general examples of
quasialgebras defined by a cochain given by a normalized Hadamard
matrix.

\begin {propos} If we consider the symmetric Hadamard matrix

$$ H=\left[\matrix{1&\hfill 1&\hfill 1&\hfill 1\cr
                     1&\hfill 1&\hfill -1&\hfill -1\cr 1&\hfill
                     -1&\hfill 1&\hfill -1\cr 1&\hfill -1&\hfill
                     -1&\hfill 1\cr}\right]
$$
as a cochain on the group $\Z_2\times \Z_2$ we obtain a nonsimple,
associative and commutative algebra.
\end {propos}
 \proof The algebra defined by $H$ and the group $(\Z_2)^2$ has the
multiplication table
$$ \matrix{e&\hfill
x&\hfill y&\hfill z\cr
        x&\hfill e&\hfill -z&\hfill -y\cr y&\hfill -z&\hfill
        e&\hfill -x\cr
          z&\hfill -y&\hfill -x&\hfill e\cr}
$$
By straightforward calculations, one can see that this algebra is
commutative and associative. Moreover, $I=<e+x, y-z>$ is an ideal
of this algebra.  \endproof

Next we consider $q=p^r$ where p is an odd prime number.  Let
$\chi(x)$ be the character defined in the finite field $F_q$ where
$\chi(0)=0$, $\chi(x)=1$ if $x$ is a square, and $\chi(x)=-1$ if
$x$ is not a square. Consider the matrix $Q$ defined by
$Q_{ij}=\chi(e_i-e_j)$ (where $e_1,e_2,e_3,...e_{q}$ are a natural
enumeration of the elements of $F_q$). Suppose that $p^r=3(mod 4)$,
that is $p^r+1=4n$ for some $n$. In this case $Q$ is skew-symmetric
and
$$
\left[\matrix{1&1&\ldots&1\cr
                -1& & & \cr \vdots&\multispan3\hfil $Q+I$\hfil\cr -1& &
                & \cr}\right]
$$
is a well-known skew symmetric Hadamard matrix  with many
applications (eg in the theory of Payley matrices). If we multiply
each row except the
first by (-1) we obtain a normalized Hadamard matrix of order $4n$,
which we can use as a cochain. In view of Proposition~3.5 we consider
for example the case where $n$ is a power of 2, so $p^r+1=2^m$ and
$m\ge 2$, and  we let
$G=(\Z_2)^m$. The octonions and quaternions appear in
this family for $p=7,3$ and $r=1$.

\begin{propos} The algebras $k_F\Z_2^m$ where  $m\ge 2$ and
$$ H=\left[\matrix{1&1&\ldots&1\cr
                1& & & \cr \vdots&\multispan3\hfil $-Q-I$\hfil\cr 1& &
                & \cr}\right]
$$
is viewed as a cochain, i.e.
\noindent $F(x,x)=-1$\quad $(x\not=e=e_0)$
\noindent$F(e_i,e_j)=-\chi(e_{i}-e_{j})\quad (i>0)\quad
(j>0)\quad (i\not=j)$, are simple and noncommutative.
\end{propos}
\proof By computing $\CR(x,y)=F(x,y)/F(y,x)$ it is easy to see that
the $k_FG$ are altercomutative (and hence noncommutative). From
this, the choice $G=(\Z_2)^m$ and $F(x,x)=-1$ for $x\not=e$, we
conclude by Proposition~3.5 that the algebras $k_FG$ admit the
diagonal strong involution $\sigma(x)=F(x,x)e$, and hence by
Proposition~3.6 that they are simple. \endproof

It is known that the Hadamard matrices have only order 2 or $4n$.
However, we can construct  another class of algebras that have
order 3 or $4n+1$ if we define the cochain by a matrix of the form

$$ \left[\matrix{1&1&\ldots&1\cr
                1& &  &\cr
           \vdots&\multispan3\hfil $H$\hfil\cr 1& &
                & \cr}\right]
$$ where $H$ is an Hadamard matrix.
The simplest case is with $H$ an Hadamard matrix of order two:

\begin{propos}
Let $k_F{\Z_3}$ be a quasialgebra which cochain defined by the
matrix,
$$ \left[\matrix{1&1&1\cr
                1 &\multispan2\hfil $H$\hfil\cr
            1 & & \cr}\right]
$$ where $H$ is an Hadamard matrix of order 2. Then $k_F{\Z_3}$ is a
nonassociative simple algebra (commutative or altercommutative).
\end{propos}

\proof
If $H$ is an Hadamard matrix of order 2 we have only two diferent
possibilities for $k_F {\Z_3}$,  with the multiplication table I:

$$ \matrix{e&\hfill
x&\hfill y\cr
        x&\hfill -y&\hfill e\cr y&\hfill e&\hfill
        x\cr}
$$ or the multiplication table II:

$$ \matrix{e&\hfill
x&\hfill y\cr
        x&\hfill y&\hfill -e\cr y&\hfill e&\hfill
        x\cr}
$$Let $I\not=k_F {\Z_3}$ be an ideal of $k_F{\Z_3}$  and let
$a=\sum_{x\in{\Z_3}} \alpha_x x$ an element of $I$. Then in both
cases $(a\cdot x)\cdot y \in I$ and if $k_F {\Z_3}$ is defined  by
table I this implies that $\alpha_x=0$ ( if not, $x\in I$ and
$I=k_F {\Z_3})$. And if $k_F {\Z_3}$ is defined by table II this
implies that $\alpha_e=0$ (if not $e\in I$ and $I=k_F {\Z_3}$). So
$I$ must be zero.
 On the other hand it is clear that the algebra defined by table I is
commutative and the one defined by the table II is
altercommutative.
\endproof

Following this idea, one also has
\begin{propos}
Let $k_F {\Z_5}$ be a quasialgebra which cochain $F$ defined by the
matrix
$$ \left[\matrix{1&1&\ldots&1\cr
                1& & & \cr \vdots&\multispan3\hfil $H$\hfil\cr 1& &
                & \cr}\right]
$$ where  $H$ is the Hadamard matrix of order 4 with maximal excess
($H=J-2I$ where $J$ has order 4 and has all elements equal to 1)
Then $k_F {\Z_5}$ is a simple commutative nonassociative algebra.
\end{propos}

One may prove this directly. Related to it is the following general
observation:

\begin{propos} Let $n\ge 4$ and $F(x,y)=(-1)^{\delta_{x,y}}$ for
$x,y\not=e$, then $k_F\Z_n$ is a simple commutative, nonassociative
algebra.
\end{propos}
\proof It's easy to see that these algebras are commutative and
nonassociative. To prove that they are simple we recall that as
$\Z_n$ is a group we know that for $x,y\in \Z_n$  we have $L_x\cdot
L_y=L_{xy}$ and so $L_x\cdot L_y\cdot L_{xy}=Id$ if we denote by
$L_x$ the map $L_x(a)=xa$. So let's consider an ideal $I\not= k_F
{\Z_n}$ and let $a=\sum_{x\in {\Z_n}}\alpha_x x$ be an element of
$I$. Now let's consider an element $y\in \Z_n$ that is diferent
from  $e$ and has order diferent from $2,3$.
 Then $y\cdot a=\sum_{x\in {\Z_n}}\alpha_x F(y,x) yx =\sum_{x\not= y\in
{\Z_n}}\alpha_x yx- \alpha_y y^2$ is an element of $I$. Now if we
multiply the last element by $z=y^2$ we will have that the element
$z\cdot (\sum_{x\not= y\in {\Z_n}}\alpha_x yx- \alpha_y
y^2)=\sum_{x\not= y\in {\Z_n}}\alpha_x y^3x+ \alpha_y y^2$ is an
element of $I$. Finally, we multiply the last element by $v=y^3$
and have that the element $v\cdot (\sum_{x\not= y \in
{\Z_n}}\alpha_x y^3x+ \alpha_y y^2)=-\alpha_e e+\sum_{x\not= e\in
{\Z_n}}\alpha_x x$ is an element of $I$.
 So $\alpha_e e\in I$ and  $\alpha_e=0$  and $I=0$ or $\alpha_e\not= 0$ and
$I=k_F {\Z_n}$.
\endproof

Other natural classes of examples (to be considered elsewhere)
include the natural generalisation of the octonions based on Galois
sequences in \cite[Chap. 2.]{Dix:div}. For example, one may take
$G=\Z_3\times\Z_3$ and $F$ of the form $F=e^{{2\pi\imath\over
3}f}$, in contrast to the $F^2=1$ case.

\section{Automorphism quasi-Hopf algebras}

We return to the general setting of Section~2, where $G$ is
equipped with a 3-cocycle $\phi$ forming a dual quasi-Hopf algebra.

We let $A$ be a finite-dimensional $G$-graded quasi-algebra in the
sense of Definition~2.3, and introduce a general construction for
its comeasuring or `automorphism' dual quasi-Hopf algebra. We let
$\{e_i \}$ be a basis of $A$ with homogeneous degrees, denoted
$|e_i |=|i|\in G$. We let $e_i \cdot e_j  =\sum_k c_{ij}{}^k e_k
$ define the structure constants of $A$ in this basis. Also, we
consider the group $G\times G$ with cocycle
\[ \phi((a,g),(b,h),(c,f))={\phi(g,h,f)\over\phi(a,b,c)}.\]

\begin{propos} Associated to a $G$-graded quasialgebra $A$ is a dual
quasi-Hopf algebra $M_1(A)$ (the comeasuring dual quasi-Hopf
algebra) defined as the free $G\times G$-quasialgebra generated by
$\{1,t^i{}_j\}$ where $|t^i{}_j|=(|i|,|j|)$ and $|1|=(e,e)$, modulo
the additional relations
\[ \sum_a c_{ij}{}^a t^k{}_a=\sum_{a,b} c_{ab}{}^k t^a{}_i\cdot t^b{}_j.\]
We define $\Delta,\eps$ as
\[ \Delta t^i{}_j=\sum_a t^i{}_a\tens t^a{}_j,\quad
\eps(t^i{}_j)=\delta^i{}_j\]
extended multiplicatively, and  extend $\phi$ to a linear
functional $\phi:M_1(A)^{\tens 3}\to k$ by
\cmath{ \phi(t^{i_1}{}_{p_1}\cdots t^{i_\alpha}{}_{p_\alpha},
t^{j_1}{}_{q_1}\cdots t^{j_\beta}{}_{q_\beta},t^{k_1}{}_{r_1}\cdots
t^{i_\gamma}{}_{p_\gamma})\\
\qquad\qquad=\delta^{i_1}{}_{p_1}\cdots
\delta^{i_\alpha}{}_{p_\alpha}\delta^{j_1}{}_{q_1}\cdots
\delta^{j_\beta}{}_{q_\beta}
\delta^{k_1}{}_{r_1}\cdots\delta^{i_\gamma}{}_{p_\gamma}
\phi(|i_1|\cdots |i_\alpha|,|j_1|\cdots |j_\beta|,|k_1|\cdots
|k_\gamma|).}
\end{propos}
\proof Since $G\times G$-graded spaces form a monoidal category,
we define the free tensor algebra on the vector space spanned by
basis $\{t^i{}_j\}$ in the usual way in a monoidal category. This
means iterated tensor products in the generators, which we
understand as nested to the right. The product is the tensor
product composed with the appropriate associativity morphism. The
$G\times G$-degree is multiplicative.  In our case the result is
the algebra generated by $1,t^i{}_j$ and the associativity rule
\cmath{\left((t^{i_1}{}_{p_1}\cdots t^{i_\alpha}{}_{p_\alpha})\cdot
(t^{j_1}{}_{q_1}\cdots t^{j_\beta}{}_{q_\beta})\right)\cdot
(t^{k_1}{}_{r_1}\cdots t^{i_\gamma}{}_{p_\gamma})\qquad\qquad\\
\qquad =(t^{i_1}{}_{p_1}\cdots
t^{i_\alpha}{}_{p_\alpha})\cdot \left((t^{j_1}{}_{q_1}\cdots
t^{j_\beta}{}_{q_\beta})\cdot (t^{k_1}{}_{r_1}\cdots
t^{i_\gamma}{}_{p_\gamma})\right) {\phi(|p_1|\cdots
|p_\alpha|,|q_1|\cdots |q_\beta|,|r_1|\cdots |r_\gamma|)
\over \phi(|i_1|\cdots
|i_\alpha|,|j_1|\cdots |j_\beta|,|k_1|\cdots |k_\gamma|)}} where
the degree of $|t^{i_1}{}_{p_1}\cdots
t^{i_\alpha}{}_{p_\alpha}|=(|i_1|\cdots |i_\alpha|,|p_1|\cdots
|p_\alpha|)$ does not depend on the nesting of the products in the
expression.

On this free quasi-associative algebra we define $\Delta,\eps$ as
shown. They are extended to products as algebra maps for the
non-associative product. It is easy to see that $\Delta,\eps$ are
compatible with the quasi-associativity, and that the extended
$\phi$ as shown makes the free quasiassociative algebra into a dual
quasi-Hopf algebra $\tilde M_1$ in the sense of Section~2. That the
extended $\phi$ is a cocycle reduces to $\phi$ a group cocycle. The
quasi-associativity axiom for a dual quasi-Hopf algebra reduces to
the $G\times G$-quasiassociativity.

Next, it is easy to verify that the quotient by the relations shown
is consistent with the $G$-quasiassociativity of our algebra $A$.
In terms of its structure constants, the latter is
\[ \sum_a c_{ij}{}^a c_{ak}{}^b=\sum_a c_{jk}{}^a
c_{ia}{}^b\phi(|i|,|j|,|k|).\] Then,
\align{(t^i{}_p\cdot t^j{}_q)\cdot
t^k{}_r\phi(|i|,|j|,|k|)c_{jk}{}^ac_{ia}{}^b
\equad&&=(t^i{}_p\cdot t^j{}_q)\cdot t^k{}_r c_{ij}{}^ac_{ak}{}^b\\
&&=c_{pq}{}^c t^a{}_c\cdot t^k{}_r c_{ak}{}^b =c_{pq}{}^c
c_{cr}{}^d t^b{}_d =c_{qr}{}^c c_{pc}{}^d\phi(p,q,r)t^b{}_d\\
&&=c_{qr}{}^c c_{pc}{}^d t^b{}_d\phi(|p|,|q|,|r|)
=c_{qr}{}^c t^i{}_p\cdot t^a{}_c c_{ia}{}^b\phi(|p|,|q|,|r|)\\
&&=t^i{}_p\cdot(t^j{}_q\cdot t^k{}_r)\phi(|p|,|q|,|r|) c_{jk}^a
c_{ia}{}^b} as required. We use the summation convention for
upper-lower indices.

Finally, we verify that these relations are compatible with
$\Delta$. Thus,
\[\Delta (c_{ab}{}^k t^a{}_i\cdot t^b{}_j)
=c_{ab}{}^k t^a{}_p \cdot t^b{}_q\tens t^p{}_i\cdot t^p{}_j
=c_{pq}{}^b t^k{}_b\tens t^p{}_i\cdot t^q{}_j=c_{ij}{}^a
t^k{}_b \tens t^b{}_a
=c_{ij}{}^a\Delta t^k{}_a\]
as required. Compatibility with $\eps$ is trivial. It is also clear
that $\phi$ restricts to the quotient in virtue of $G$-grading of
the algebra $A$. For example, $\phi(t^i{}_j,t^a{}_p t^b{}_q
c_{ab}{}^k,t^r{}_s)=\delta^i{}_j\delta^r{}_s\phi(|i|,|p||q|,|r|)
c_{pq}{}^k$ while $\phi(t^i{}_j c_{pq}{}^a
t^k{}_a,t^r{}_s)=\delta^i{}_j\delta^r{}_s\phi(|i|,|k|,|r|)
c_{pq}{}^k$, but $c_{pq}{}^k=0$ unless $|p||q|=|r|$. Thus, the
quotient of $\tilde M_1$ by the relations shown defines a dual
quasi-Hopf algebra $M_1$. \endproof

If $G$ is in addition equipped with a quasi-bicharacter $\CR$
(making $(kG,\phi)$ into a dual quasitriangular quasi-Hopf algebra)
then there is a natural braiding in the category of $G$-graded
spaces as explained in Section~2.  In our case, we extend $\CR$ to
a quasi-bicharacter on $G\times G$ by
\[ \CR((a,g),(b,h))={\CR(g,h)\over\CR(a,b)}\]

\begin{propos} If $\CR$ is a quasi-bicharacter on $G$, the
comeasuring dual quasi-bialgebra $M_1(A)$ has a natural quotient
$M_1(\CR,A)$ with the additional relation of quasi-commutativity as
a $G\times G$-quasialgebra. Then $M_1(\CR,A)$ is a
dual-quasitriangular dual quasi-Hopf algebra with $\CR$ extended as
a linear functional on $M_1(\CR,A)^{\tens 2}$ by
\[\CR(t^{i_1}{}_{p_1}\cdots t^{i_\alpha}{}_{p_\alpha},
t^{j_1}{}_{q_1}\cdots t^{j_\beta}{}_{q_\beta})
=\delta^{i_1}{}_{p_1}\cdots
\delta^{i_\alpha}{}_{p_\alpha}\delta^{j_1}{}_{q_1}\cdots
\delta^{j_\beta}{}_{q_\beta}
\CR(|i_1|\cdots |i_\alpha|,|j_1|\cdots |j_\beta|) \]
\end{propos}
\proof Explicitly, quasicommutativity as a $G\times G$-graded algebra is
\cmath{ (
t^{j_1}{}_{q_1}\cdots t^{j_\beta}{}_{q_\beta})\cdot
(t^{i_1}{}_{p_1}\cdots t^{i_\alpha}{}_{p_\alpha})\CR(|p_1|\cdots
|p_\alpha|,|q_1|\cdots |q_\beta|)\qquad\\
\qquad=\CR(|i_1|\cdots
|i_\alpha|,|j_1|\cdots |j_\beta|) (t^{i_1}{}_{p_1}\cdots
t^{i_\alpha}{}_{p_\alpha})\cdot ( t^{j_1}{}_{q_1}\cdots
t^{j_\beta}{}_{q_\beta})}
 for the braiding $\Psi$ determined by
$\CR$ on $G\times G$ as shown. This is the quasicommutativity
property of a dual-quasitriangular dual quasi-Hopf algebra in the
sense of Section~2, with $\CR$ defined as stated. That this $\CR$
is well-defined on the free quasi-associative algebra $\tilde M_1$
is clear. That it descends to the quotient by the relations of
$M_1$ follows by the $G$-grading of our algebra $A$ as for $\phi$
in the preceding proof. That it well-defined on $M_1(\CR,A)$ itself
requires repeated use of the quasi-bicharacter property and is
omitted
\note{***}\endproof

Moreover, both $M_1(A)$ and hence $M_1(\CR,A)$ coact on $A$:

\begin{propos} $M_1$ coacts on $A$ by $\beta:e_i \mapsto
e_a\tens t^a{}_i$ and $\beta$ is an algebra map.
\end{propos}
\proof The definition of the coaction is consistent with the
relations of $A$:
\[ \beta(e_i \cdot e_j  )=\beta(e_i )\beta(e_j  )=e_a\cdot e_b\tens
t^a{}_i\cdot t^b{}_j
=c_{ab}{}^k e_k    \tens t^a{}_i\cdot t^b{}_j=c_{ij}{}^a e_k
\tens t^k{}_a
=\beta(c_{ij}{}^ae_a)\]
in virtue of the relations of $M_1$. \endproof

In fact, it should be clear from the proof that the relations of
$M$ are the minimum relations such that a coaction of this form
extends as an algebra map. In the case where $\phi$ is trivial we
recover in fact the dual (arrows-reversed) version of the measuring
bialgebra $M(A,A)$ in \cite{Swe:hop}, and this is the reason
motivation behind our construction. Some further recent
applications of this comeasuring bialgebra construction in the
associative (not quasi-associative) case appear in \cite{Ma:aut}.

Before turning to examples, we note that $M_1(A)$ and $M_1(\CR,A)$
have natural further quotients. Thus

\begin{propos} The diagonal quotient $M_D$ of $M_1$ by the relations
$t^i{}_j=0$ when $i\ne j$ is an associative Hopf algebra with
generators $t_i$ and relations
\[ c_{ij}{}^k(t_k-t_it_j)=0,\quad \Delta(t_i)=t_i\tens t_i,\quad
\eps(t_i)=1.\]
It can also be viewed as a dual quasi-Hopf algebra with
\[ \phi(t_{i_1}\cdots t_{i_\alpha},t_{j_1}\cdots
t_{j_\beta},t_{k_1}\cdots t_{k_\gamma})=\phi(|i_1|\cdots
|i_\alpha|,|j_1|\cdots |j_\beta|,|k_1|\cdots |k_\gamma|)\] The
diagonal quotient $M_D(\CR,A)$ of $M_1(A)$ is the commutative
quotient of $M_D$, and can be viewed as a dual quasitriangular dual
quasi-Hopf algebra with \[ \CR(t_{i_1}\cdots
t_{i_\alpha},t_{j_1}\cdots t_{j_\beta})=\CR(|i_1|\cdots
|i_\alpha|,|j_1|\cdots |j_\beta|).\]
\end{propos}
\proof This is elementary. The quasi-associativity and quasicommutativity
of $M_1(A),M_1(\CR,A)$ clearly reduce in the diagonal case to usual
associativity and commutativity. The coproduct becomes group-like.
\endproof

Finally, we have not required yet that $A$ is unital. When it is,
we choose our basis so that $e_0=1$. In this case we let $\{e_i \}$
denote the remaining basis elements.

\begin{propos} When $e_0=1$ the unit of $A$, we define $M_0(A)$ as the
quotient of $M_1(A)$ by $t^0{}_0=1$ and $t^i{}_0=t^0{}_i=0$. This
forms a dual quasi-Hopf algebra with relations and coproduct
\[c_{ij}{}^at^k{}_a=c_{ab}{}^k t^a{}_i t^b{}_j,\quad c_{ij}{}^0c
=c_{ab}{}^0 t^a{}_i t^b{}_j,\quad \Delta t^i{}_j=t^i{}_a\tens t^a{}_j.\]
Similarly, $M_0(\CR,A)$ remains dual quasitriangular and preserves
this form.
\end{propos}
\proof We set $t^i{}_0=0=t^0{}_i$ and denote $t^0{}_0=c$. Note that
 $c_{i0}{}^j=c_{0i}{}^j=\delta^j{}_i$
and $c_{00}{}^0=1$. Therefore, the relations of $M_1(A)$ become
(for all labels not $0$), the relations as stated, and the
additional relations
\[ t^i{}_j c=t^i{}_j=ct^i{}_j,\quad c^2=c.\]
In view of these latter relations, it is natural to set $c=1$.
Moreover, the matrix coproduct and counit are clearly compatible
with the quotient. One has $\Delta c=c\tens c$, so this is
consistent with $c=1$ as well. Moreover, since $|e_0|=e$, the
identity in $G$, it is clear that setting $c=1$ is consistent with
the definition of $\phi,\CR$ on $M_1$. Likewise, their
delta-function form is consistent with $t^i{}_0=0=t^0{}_i$. Hence
$M_0(A)$ and $M_0(\CR,A)$ inherit these structures and are dual
quasi-Hopf algebras.

Finally, the coaction of $M_1$ becomes $\beta(e_0)=e_0\tens c$ and
$\beta(e_i )=e_a\tens t^a{}_i$. The relations of $M_0$ are such
that the bilinear form on the span of $\{e_i \}$  defined by
$c_{ij}{}^0$ is preserved.
\endproof

Given our basis, we can identify $A/1$ with the span of $\{e_i \}$
for $i\ne 0$, and $B(e_i ,e_j  )=c_{ij}{}^0$ is a natural bilinear
form on it. We see that our reduced comeasuring dual quasi-quantum
groups $M_0(A)$, $M_0(\CR,A)$ preserve this. Also, the two
relations for the $t^i{}_j$ imply that
\[ c_{ab}{}^dc_{dc}{}^0 (t^a{}_i\cdot t^b{}_j)\cdot t^c{}_k
=c_{dc}{}^0c_{ij}{}^a t^d{}_a\cdot t^c{}_k
=c_{ij}{}^a c_{ak}{}^0\]
so that the trilinear form $c_{ij}{}^a c_{ak}{}^0$ is also
preserved in a certain sense. Finally, we have the further diagonal
quotients of $M_0(A)$ and $M_0(\CR,A)$.

\begin{corol} For $F$ a cochain on $G$, and basis $G$ of $A=k_FG$, the
diagonal quotient $M_{D0}=(kG,\phi)$ as a dual quasi-Hopf algebra.
\end{corol}
\proof In this basis $c_{ij}{}^k=1$ iff $ij=k$ in $G$ and zero otherwise.
Hence the relations of $M_D$ are $t_it_j=t_k$ for $k=ij$ and empty
for $k\ne ij$. In $M_{D0}$ we further identify $t_e=1$ as in the
group algebra. Finally, $|t_i|=i$ and we obtain $kG,\phi$ as a dual
quasi-Hopf algebra. When $G$ is commutative we obtain a commutative
algebra and $M_{D0}=M_{D0}(\CR,k_FG)$. \endproof

We now compute these constructions for the complex numbers and for
the quaternions, as real two and four dimensional algebras. More
generally, we work over a general ground field of characteristic
not 2.

\begin{example} When $A=k[i]$ (where $i^2=-1$), the comeasuring
bialgebra $M_1(k[i])$ is generated by $1$ and a matrix of
generators $\pmatrix{a&b\cr c&d}$, with the relations
\[ a^2-c^2=a=d^2-b^2,\quad ac+ca=c, \quad -c=bd+db,\quad
ab-cd=b=ba-dc,\quad ad+cb=d=bc+da.\] This has a natural bialgebra
quotient of the form $\pmatrix{c&s\cr
-s&c}$ with
\[ \Delta c=c\tens c-s\tens s,\quad \Delta s=s\tens c+c\tens s,\quad
c^2-s^2=c,\quad sc+cs=s.\]

The quotient $M_0=M_0(\CR)=M_{D0}$ is $k\Z_2$ (as generated by
$d$), and its coaction is $\beta(1)=1\tens 1,\beta(i)=i\tens d$.
\end{example}
\proof We write out the 8 relations for $M_0$ using the structure
constants of $k[i]$. The quotient $M_0$ is already diagonal and
commutative. Hence by the preceding corollary, it gives $kG=k\Z_2$.
Its coaction is on $k[i]$ is the canonical nontrivial one
corresponding to the $G$-grading. Note that evaluating with the
nontrivial character of $\Z_2$ gives the canonical automorphism
$i\to -i$. \endproof

The intermediate quotient here is the `trigonometric bialgebra':
the coproduct has the same form as the addition rules for the sine
and cosine functions. Whereas it is usually considered as a
coalgebra\cite{Swe:hop}, we obtain here a natural algebra structure
forming a bialgebra. It too coacts on $k[i]$ by our constructions
as the push out of the universal coaction of $M_1$.

We also note that when $A=k \Z_2$, the comeasuring bialgebra
$M_1(k\Z_2)$ has the same form as in the preceding example but with
all minus signs replaces by $+$. The quotient of the form
$\pmatrix{a&b\cr b&a}$ can then be diagonalised as $g^\pm=a\pm b$
and becomes the bialgebra
\[ g^\pm g^\pm=g^\pm,\quad \Delta g^\pm=g^\pm\tens g^\pm,\quad\eps g^\pm=1\]
of two mutually noncommuting projectors $g^\pm$. This is an
infinite-dimensional algebra with every element of the form either
$g^+g^-g^+\cdots$ or $g^-g^+g^-\cdots$ (alternating). One may make
a similar diagonalisation $g^\pm=c\pm i s$ for the trigonometric
bialgebra in the case when $i=\sqrt{-1}\in k$.

\begin{propos} When $A=\H$ the quaternion algebra over $k$, the comeasuring
bialgebra $M_0(\H)$ has generators $1$ and three vectors of
generators $\vec{t_j}=(t^i{}_j)$, $i=1,2,3$ and relations
\[ \vec{t_1}\times\vec{t_2}=\vec{t_3},\quad \vec{t_1}\times\vec{t_1}=0,
\quad{\rm + cyclic},\quad \vec{t_i}\cdot\vec{t_j}=\delta_{ij}. \]
Here $\times$ is the vector cross product and $\cdot$ is the vector
dot product. The quotient $M_0(\CR)$ is defined by the additional
relation that the generators commute.
\end{propos}
\proof We choose the standard basis (where $e_0=1$ and $e_i $, $i=1,2,3$
have the relations $ee_2=e_3$ and $e^2=-1$ and their cyclic
permutations). In this case the structure constants are
$c_{ij}{}^k=\eps_{ijk}$, the totally antisymmetric tensor with
$\eps_{123}=1$, and $c_{ij}{}^0=-\delta_{ij}$, the standard
Euclidean metric. The relations of $M_0$ then become
\[ \delta_{ab}t^a{}_i t^b{}_j=\delta_{ij},\quad t^a{}_i t^b{}_j\eps_{abc}
=t^c{}_a\eps_{ij}{}^a\]
plus cyclic permutations.  For $M_0(\CR)$ it is enough to note that
$\CR(|i|,|j|)=-1$ hence that the $t^i{}_j$ mutually commute. In
this case we have $\det(m)=1$, i.e. $M_0(\CR)$ is a quotient of
$k[SL_3]$.
\endproof

The relations in $M_0$ here are asymmetric, a reflection of their
role as coacting from the right on $\H$. If we consider also the
left-handed versions of our constructions, we have a joint quotient
where
\[ t^i{}_at^j{}_b \delta^{ab}=\delta^{ij}\]
is added. In this case we have a Hopf algebra with
$St^i{}_j=t^j{}_i$ and the corresponding quotient of $M_0(\CR)$ is
a quotient of the group coordinate ring $k[SO_3]$. Since $M_0(\CR)$
is universal among commutative bialgebras coacting on $\H$, it
follows that, it must project onto the group coordinate ring of the
classical automorphism group.

The corresponding computation for the octonions yields for $M_0$ a
dual quasi-Hopf algebra with nontrivial $\phi$. Its detailed form
is somewhat more complex than the quaternion case, however; on
general grounds we know that it projects for example on to the
group coordinate ring $k[G_2]$ (the classical automorphism Hopf
algebra of the octonions).
\note{***}

\section{Quasiassociative linear algebra}

In this section we use our categorical approach to octonions to
provide the natural `quasiassociative' setting for the basic linear
algebra associated to them. We define the natural notion of
`representation'. We also provide the definition of $V^*$ for any
finite-dimensional $G$-graded vector space, and the associated
endomorphism quasialgebra $V\tens V^*$. These constructions are the
specialization to the $G$-graded quasi-algebra setting of standard
constructions for braided categories.

Thus, the notion of representations, indeed of all linear algebra
and quantum group constructions, make sense in any braided
category, see \cite{Ma:introm}. One writes all constructions as
compositions of morphisms, inserting the associator $\Phi$ as
necessary. For example, in the case of the category of
$(kG,\phi)$-comodules, we clearly have:

\begin{defin} A representation or `action' of a $G$-graded
quasialgebra $A$ is a $G$-graded vector space $V$ and a
degree-preserving map $\la:A\tens V\to V$ such that
\[ (ab)\la v=\phi(|a|,|b|,|v|)a\la(b\la v),\quad 1\la v=v\]
on elements of homogeneous degree. Here $|a\la v|=|a||v|$.
\end{defin}

This is the obvious polarization of the quasi-associativity of the
product of $A$. Clearly, a quasialgebra acts on itself by the
product map (the regular representation).

Next, we recall that an object $V$ in a braided category is called
`rigid' if there is an object $V^*$ and morphisms
\[ \ev:V^*\tens V\to \und 1,\quad \coev:\und 1\to V\tens V^*\]
such that
\[ (\id\tens\ev)\Phi_{V,V^*,V}(\pi\tens\id)=\id,\quad (\ev\tens\id)
\Phi^{-1}_{V^*,V,V^*}(\id\tens\pi)=\id \]
holds. In the case of the comodule category of a dual quasiHopf
algebra, these maps exist whenever $V$ is finite-dimensional, see
\cite{Ma:book} for the explicit formulae, cf\cite{Dri:qua}. For the
dual quasiHopf algebra $kG,\phi$, i.e. for the category of
$G$-graded vector spaces, these maps are given by
\[ \ev(f^i\tens e_j  )=\delta^i{}_j,\quad \coev(1)=\sum_i e_i \tens
f^i\phi^{-1}(|i|,|i|^{-1}, |i|)\] in terms of a basis of $V$ with
degree $|e_i |\equiv |i|\in G$ and its usual dual, i.e. $\ev$ can
be taken as the usual evaluation, $V^*$ as the usual dual, but
$\coev$ is modified by the group 3-cocycle. Here $|f^i|=|i|^{-1}$
so that $\ev,\coev$ are degree preserving.

\begin{propos}
If $V$ is rigid then $\End(V)=V\tens V^*$ becomes a $G$-graded
quasialgebra. The product map is
\[ (v\tens f)(w\tens h)=v\tens h <f,w> {\phi(|v|,|f|,|w||h|)\over
\phi(|f|,|w|,|h|)}.\]
Moreover, finite-dimensional representations $V$ of a quasialgebra
$A$ are in 1-1 correspondence with quasialgebra maps $A\to
\End(V)$.
\end{propos}
\proof This is proven by commuting diagrams exactly as one would
prove these statements in linear algebra, only inserting the
associator $\Phi$ wherever needed to change bracketing. Thus (in
any monoidal category) one finds
\[ (\id\tens(\ev\tens\id))\circ(\id\tens \Phi^{-1}_{V^*,V,V^*})
\circ\Phi_{V,V*,V\tens V^*}
:(V\tens V^*)\tens (V\tens V^*)\to V\tens V^*\]
as the natural product on $\End(V)=V\tens V^*$. Its action on $V$
is the map
\[(\id\tens\ev)\circ\Phi_{V,V^*,V}:(V\tens V^*)\tens V\to V.\]
If $A$ acts on $V$ then
\[ \rho=(\la\tens\id)\circ\Phi_{A,V,V^*}\circ(\id\tens\coev):
A\to V\tens V^*\]
is an algebra map. Conversely, given $\rho$,
\[ (\id\tens\ev)\circ\Phi_{V,V^*,V}\circ(\rho\tens\id):A\tens V\to V\]
is an action on $V$. We then specialize to the case of $G$-graded
quasialgebras using the form of $\Phi$ in terms of the 3-cocycle
$\phi$. \endproof

In lieu of all the commutative diagrams here, we will prove this
more explicitly in a concrete form for our particular setting.
First, we identify $V\tens V^*$ with matrices in the usual way
relative to our basis, i.e.
\[ \alpha=\sum \alpha^i{}_j E_i{}^j,\quad E_i{}^j=e_i\tens f^j\]
as a definition of the components of $\alpha\in V\tens V^*$. Then
the preceding proposition translates into the following
proposition. We write $n=\dim(V)$ and $|i|\in G$ as the further
data provided by $V$, which we also use.

\begin{propos} Let $|i|\in G$ for $i=1,\cdots,n$ be a choice of
grading function. Then the usual $n\times n$ matrices $M_n$ with
the new product
\[ (\alpha\cdot \beta)^i{}_j=\sum_k \alpha^i{}_k
\beta^k{}_j{\phi(|i|,|k|^{-1},|k||j|^{-1})\over
\phi(|k|^{-1},|k|,|j|^{-1})},\quad \forall\alpha,\beta\in M_n\]
form a $G$-graded quasialgebra $M_{n,\phi}$, where
$|E_i{}^j|=|i||j|^{-1}\in G$ is the degree of the usual basis
element of $M_n$. An action of a $G$-graded quasialgebra in the
$n$-dimensional vector space with grading $|i|$ is equivalent to an
algebra map $\rho:A\to M_{n,\phi}$.
\end{propos}
\proof The product suggested by the preceding proposition
is
\[ E_i{}^j\cdot E_k{}^l=\delta^j_kE_i{}^l{\phi(|i|,|j|^{-1},
|j||l|^{-1})\over
\phi(|j|^{-1},|j|,|l|^{-1})},\]
which yields the formula shown for $\alpha=\sum\alpha^i{}_j
E_i{}^j$, etc. This product is quasiassociative since
\align{\equad && (E_i{}^j\cdot E_k{}^l)\cdot E_m{}^n=
\delta^j_k\delta^l_m E_i{}^n{\phi(|i|,|m|^{-1},|m||n|^{-1})
\phi(|i|,|k|^{-1},|k||l|^{-1})\over\phi(|m|^{-1},|m|,|n|^{-1})
\phi(|k|^{-1},|k|,|l|^{-1})}\\
\equad&& E_i{}^j\cdot(E_k{}^l\cdot E_m{}^n)\phi(|E_i{}^j|,
|E_k{}^l|,|E_m{}^n|)\qquad\\
&&=\delta^l{}_m\delta^j{}_k{\phi(|i||j|^{-1},|k||l|^{-1},|m||n|^{-1})
\phi(|k|,|m|^{-1},|m||n|^{-1})\phi(|i|,|k|^{-1},|k||n|^{-1})\over
\phi(|m|^{-1},|m|,|n|^{-1})\phi(|k|^{-1},|k|,|n|^{-1})}}
which are equal since
\[ \phi(|i||k|^{-1},|k||m|^{-1},|m||n|^{-1})={\phi(|k|^{-1},|k|,
|n|^{-1})\phi(|i|,|k|^{-1},|k||m|^{-1})
\phi(|i|,|m|^{-1},|m||n|^{-1})\over
\phi(|k|^{-1},|k|,|m|^{-1})\phi(|k|,|m|^{-1},|m||n|^{-1})
\phi(|i|,|k|^{-1},|k||n|^{-1})}\]
by repeated use of the cocycle property of $\phi$.

Note also that the grading function $|i|$ is equivalent to
specifying a $G$-graded vector space $V=\{e_i\}$ with grading
$|e_i|=|i|$. An action of a $G$-graded quasialgebra $G$ is
equivalent to structure constants $v_{\alpha i}{}^j$ such that
\[  c_{\alpha\beta}{}^\gamma v_{\gamma i}{}^j
=v_{\beta i}{}^k v_{\alpha k}{}^j\phi(|\alpha|,|\beta|,|i|),
\quad v_{0i}{}^j=\delta^j_i\]
where  $\{x_\alpha\}$ (say) is a basis of $A$ with $x_0=1$ and
$x_\alpha\la e_i=v_{\alpha i}{}^j e_j$. The corresponding map $A\to
M_{n,\phi}$ is
\[ \rho(x_\alpha)^i{}_j={v_{\alpha j}{}^i\over \phi(|j|,
|j|^{-1},|j|)\phi(|i||j|^{-1}, |j|,|j|^{-1})}.\] That $\rho$ is an
algebra map is (from the definitions stated) is
\[ c_{\alpha\beta}{}^\gamma v_{\gamma j}{}^i {1\over
\phi(|j|,|j|^{-1},|j|)
\phi(|i||j|^{-1},|j|,|j|^{-1})}\qquad\qquad\qquad\]
\[=v_{\alpha a}{}^i v_{\beta j}{}^a {\phi(|i|,|a|^{-1},|a||j|^{-1})
\over \phi(|a|^{-1},|a|,|j|^{-1})\phi(|a|,|a|^{-1},|a|)
\phi(|i||a|^{-1},|a|,|a|^{-1})
\phi(|j|,|j|^{-1},|j|)\phi(|a||j|^{-1},|j|,|j|^{-1})}.\]
Since the structure maps are degree preserving, we know that
$|a|=|\beta||j|,\quad |i|=|\alpha||a|$ for nonzero terms on the
right and side. That $\rho$ is an algebra map is then equivalent to
$v_{\alpha i}{}^j$ an action in view of the identity
\align{&&\equad \phi(|\alpha||\beta|,|j|,|j|^{-1})
\phi(|\alpha||\beta||j|,|j|^{-1}|\beta|^{-1},|\beta|)\\
&&=
\phi(|\beta|,|j|,|j|^{-1})\phi(|\beta||j|,|\beta|^{-1}
|j|^{-1},|\beta||j|)
\phi(|\alpha|,|\beta||j|,|\beta|^{-1}|j|^{-1})
\phi(|\beta|^{-1}|j|^{-1},|\beta||j|,|j|^{-1})
\phi(|\alpha|,|\beta|,|j|)}
which holds by repeated use of the 3-cocycle property. \endproof

For example, the left regular representation of a quasialgebra on
itself provides  representation $\rho:A\to M_{n,\phi}$ where
$n=\dim(A)$. For the octonions, for example, we have a
representation in $8\times 8$ quasimatrices.


\end{document}